\documentclass{amsart}
\usepackage[utf8]{inputenc}
\usepackage{amsmath}
\usepackage{amsthm}

\newtheorem*{theorem}{Theorem}

\newcommand\id{\mathrm{id}}

\title{Cyclic extensions are radical}
\date{November 8, 2014}

\author{Mariano Suárez-Álvarez}
\email{mariano@dm.uba.ar}
\address{Departamento de Matem\'atica.
  Facultad de Ciencias Exactas y Naturales.
  Universidad de Buenos Aires.
  Ciudad Universitaria, Pabell\'on I.
  (1428) Ciudad de Buenos Aires (Argentina)}
\thanks{This work was partially supported by  UBACYT X475, PIP-CONICET
  2012-2014 11220110100870, PICT 2011-1510 and MathAmSud-GR2HOPF.}


\begin{document}

\maketitle

\thispagestyle{empty}

The fact that finite Galois extensions with cyclic Galois group can be
constructed by adjoining a root of a well-chosen element of the base field is
a basic result of Galois theory and a key point in studying the problem of
solvability by radicals. In textbooks on the subject, this result is
usually obtained as a consequence of Hilbert's Theorem~90, which is itself
deduced from Artin's theorem on the independence of characters. The very
simple argument we present below sidesteps those requirements.

\begin{theorem}
If $E/K$ is a cyclic extension of degree~$n$ of a field~$K$ which contains
a primitive $n$th root of unity~$\zeta$, then there exists an  $x\in
E$ with $x^n\in K$ and $E=K(x)$.
\end{theorem}

Notice that the existence of a primitive $n$th root of unity in~$K$ implies
that the characteristic of this field does not divide~$n$.

\begin{proof}
Let $\sigma$ be a generator of the Galois group~$G$ of the extension.
Since $\sigma^n=\id_E$, the minimal polynomial of~$\sigma$ over~$K$ divides
$X^n-1$ in~$K[X]$ and then, as this polynomial has all its roots in~$K$ and
they are all simple, $\sigma$ is diagonalizable and its eigenvalues are
$n$th roots of unity. 

Let $\Gamma$ be the set of eigenvalues of~$\sigma$. If
$\lambda$,~$\mu\in\Gamma$, so that there are $a$,~$b\in E^\times$ such that
$\sigma(a)=\lambda a$ and $\sigma(b)=\mu b$, then $\lambda\mu\in\Gamma$, as
$\sigma(ab)=\lambda\mu ab$. Since $\Gamma$ is contained in~$\Omega_n$, the
finite group of $n$th roots of unity, this is enough to conclude that
$\Gamma$ is in fact a \emph{subgroup} of~$\Omega_n$. If $m$ is the order
of~$\Gamma$ , then we have that $\lambda^m=1$ for all~$\lambda\in\Gamma$
and, as $\sigma$ is diagonalizable with eigenvalues in~$\Gamma$,
that $\sigma^m=\id_E$. As the order of~$\sigma$ is~$n$ and $m\leq n$, it
follows from this that $m=n$.

All $n$th roots of unity are therefore eigenvalues of~$\sigma$ and, in
particular, there is an $x\in E^\times$ such that $\sigma(x)=\zeta x$.
Then $\sigma(x^n)=\zeta^nx^n=x^n$, so $x^n$ is in the fixed field
$E^G=K$. We thus see that $x$ is a root of the polynomial
$X^n-x^n$ of~$K[X]$, which is irreducible over~$K$: indeed, it has simple
roots $x$,~$\zeta x$,~$\dots$,~$\zeta^{n-1}x$ in~$E$ and these are permuted
transitively by~$\sigma$. The degree of the subextension $K(x)/K$
of~$E/K$ is then~$n$ and, of course, this means that $E=K(x)$.
\end{proof}

\vfill

\end{document}